\documentclass[12pt]{article}
\usepackage{amssymb,amsmath}

\def\tr {{\rm tr}}

\def\dim {{\rm dim}}

\begin{document}

\title{On the space of almost complex structures}
\author{Daurtseva N.~A., Smolentsev N.~K.
\thanks{Kemerovo State University, Kemerovo, 650043, RUSSIA.
Email smolen@kemsu.ru, }
} \maketitle

\begin{abstract}
It is shown in the article \cite{Sm} that the space of associated metrics
and almost complex structures on the symplectic manifold is infinitely
dimensional Kahler manifold. The results of
\cite{Sm} are generalized on the space of all almost complex structures
of manifold.
\end{abstract}

\vspace {5mm}

{\bf 1. The space of almost complex structures.} Let $M$
is smooth closed oriented manifold of dimension $2n$.

{\it Almost complex structure} on $M$ is field of endomorphisms
$J_x:T_xM\longrightarrow T_xM$, varying smoothly with $x\in M$
satisfying the anti-involution property: $J^2=-Id$, where $Id$ is
identical endomorphism.

Let $\mathcal{A}$ is the space of all smooth almost complex
structures on $M$. It is the space of smooth sections of the
bundle $A(M)$ over $M$, with automorfisms $J_x$ of $T_xM$, such
that $J_x^2= -I_x$ as the fiber $A_x(M)$ over the point $x\in M$.
Since the space $\mathcal{A}$ is space of sections $\mathcal{A} =
\Gamma (A(M))$, then it is infinitely dimensional smooth
ILH-manifold.

Let $J\in \mathcal{A}$ is almost complex structure on $M$.
Differentiating the property $J^2= -1$, one can obtain that space
$T_J\mathcal{A}$ consists of endomorphisms $K:\ TM \to TM$,
anticommutating with $J$, $JK+KJ=0$. Denote space of such
operators on $TM$ as $\mathrm{End}_{J}(TM)$. Thus,
$$
T_J\mathcal{A}=\mathrm{End}_{J}(TM).
$$

As the space $\mathcal{A}$ is smooth ILH-manifold, one can define
local maps in usual way. There is more natural parameterization by
using the Cayley transformation.

Let $J_0$ is any fixed almost complex structure. The tangent space
$T_{J_0}\mathcal{A}$ consists of endomorphisms $K:TM\to TM$,
anticommutating with almost complex structure $J_0$,\ $K J_0 =
-J_0 K$. That is why we have $J_0e^K=e^{-K}J_0$ for the
exponential mapping $e^K$. It immediately follows, that operator
$$
J=J_0e^K
$$
is almost complex structure. The last equality gives
parameterization of the space $\mathcal{A}$ in the neiborhood of
the element $J_0$ by endomorphisms $K$, anticommutating with
$J_0$:
$$
E:\ \mathrm{End}_{J_0}(TM) \longrightarrow \mathcal{A},\quad
K\mapsto J = J_0e^K.
$$

In the matrix theory the rational dependence
$w=\frac {1+iz}{1-iz}$, $z=i\frac {1-w}{1+w}$ is used instead of
transcendental one $w=e^{iz}$.
Apply this transformation to the operator
$K$, which has property $K J_0 = -J_0 K$, we obtain
$$
J=J_0\left(1+K\right)\left(1-K\right)^{-1}.
$$
It is easy to see, that
$$
J=J_0\left(1+K\right)\left(1-K\right)^{-1}=\left(1-K\right)\left(1+K\right)^{-1}J_0.
$$
That is why $J$ is almost complex structure. When we defined such
almost complex structure we supposed  that $1-K$ is degenerated.
It is enough to this, that the operator $K(x)$ has to have no
properly numbers equaled to unit at each point $x\in M$. Clearly,
that set of such endomorphisms is open in the space
$\mathrm{End}_{J_0}(TM)$ of all endomorphisms $K:TM\to TM$,
anticommutating with $J_0$. Denote this set by $\mathcal{V}(J_0)$,
$$
\mathcal{V}(J_0)=\{K\in \mathrm{End}(TM);\ \ K J_0 = -J_0 K,\ 1-K
\ - \mbox{
invertible } \}.
$$

{\bf Proposition 1.}\cite{Sm} {\it Relations
$$
J=J_0\left(1+K\right)\left(1-K\right)^{-1}, \eqno{(1)}
$$
$$
K=\left(1-JJ_0\right)^{-1}\left(1+JJ_0\right)
$$
give one-to-one correspondence between the set of endomorphisms
$K:TM\to TM$,\\ anticommutating with almost complex structure $J_0$
such that $1-K$ is invertible, and the set of almost complex
structure $J$ on $M$ for which endomorphism $1-JJ_0$ is
invertible.}

\vspace{3mm}

{\bf Remark.} Algebraic clearer relation
$$
J=\left(1-K\right)J_0\left(1-K\right)^{-1}
$$
follows from the (1).

\vspace{3mm}

Set
$$
\mathcal{U}(J_0)=\{J\in \mathcal{A};\  1 - JJ_0 - \mbox{
isomorphism } TM \}
$$
is open in $\mathcal{A}$. That is why mapping
$$
\Phi:\ \mathcal{U}(J_0) \longrightarrow \mathcal{V}(J_0), \qquad J
\mapsto K,
$$
$$
K=(1 - JJ_0)^{-1}(1 + JJ_0) \eqno{(2)}
$$
gives local coordinates in the neiborhood of $J_0$.
If $K=\Phi(J)$, then
$$
J = J_0(1 + K)(1 - K)^{-1}.
$$
We will denote such almost complex structure as $J_K$.

\vspace{3mm}

The formulas of "changing coordinates" are easy obtained from (2).
If $J\in \mathcal{U}(J_0)\cap \mathcal{U}(J_1)$ and
$K=(1 - JJ_0)^{-1}(1 + JJ_0)$,\ $P=(1-JJ_1)^{-1}(1+JJ_1)$, then
$$
P = (1-(1-K)(1+K)^{-1}J_0J_1)^{-1}(1+(1-K)(1+K)^{-1}J_0J_1).
$$

\vspace{3mm}

Let $A\in \mathrm{End}_{J_0}(TM)$. Find element $A^*\in
T_{J_K}\mathcal{A}$, which corresponds to element $A$ under
differential of coordinate mapping:
$$
d\ \Phi^{-1}_{K}:\ T_K\mathrm{End}_{J_0}(TM) \longrightarrow
T_{J_K}\mathcal{A}.
$$
It is enough to take the curve $K_t=K+tA$ and differentiate
equality
$J_{K_t}=J_0(1+K_t)(1-K_t)^{-1}$ with respect to $t$.
$$
A^*=\left.\frac{d}{dt}\right|_{t=0}
J(K_t)=J_0\left.\frac{d}{dt}\right|_{t=0}K_t(1-K)^{-1}+
J_0(1+K)\left.\frac{d}{dt}\right|_{t=0}(1-K_t)^{-1}=
$$
$$
=J_0A(1-K)^{-1}+J_0(1+K)(1-K)^{-1}A(1-K)^{-1}=J_0(1+(1+K)(1-K)^{-1})A(1-K)^{-1}=
$$
$$
=J_0(1-K+1+K)(1-K)^{-1}A(1-K)^{-1}=2J_0(1-K)^{-1}A(1-K)^{-1}.
$$

Therefore, we obtain
$$
A^*=2J_0(1-K)^{-1}A(1-K)^{-1}.
$$
The following expressions are easy obtained:
$$
A^*=J_0\ A(1-K)^{-1} + J\ A(1-K)^{-1},
$$
$$
A^*=2J_K(1-K)(1-K^2)^{-1}A(1-K)^{-1}.
$$

\vspace{3mm}

{\bf 2. Pseudo-Riemannian structure on $\mathcal{A}$.} Fix Riemannian
structure $g_0$ on the manifold $M$. Then one can define the following
structures on the space $\mathcal{A}$:

1. The weak Pseudo-Riemannian structure. If $A,B\in T_J\mathcal{A}$,
then their inner product is defined by formula:
$$
(A,B)_{J}=\int_M \tr (A\circ B)d\mu(g_0),
$$
where $\mu(g_0)$ is Riemannian volume element.

2. Almost complex structure.
$$
{\bf J}_J:\ T_J\mathcal{A} \longrightarrow T_J\mathcal{A},\quad
{\bf J}_J(A)=A\circ J.
$$

3. Antisymmetric non-degenerated 2-form $\Omega$ on
$\mathcal{A}$. If $A,B\in T_J\mathcal{A}$, then
$$
\Omega_J(A,B)=\int_M \tr(AJB)d\mu(g_0)=({\bf J}_J A,B)_{J}.
$$

The expressions of these structures in local "coordinates"
on $\mathcal{A}$ in the neiborhood $\mathcal{U}(J_0)$ of the element $J_0$ are
easy found:
$$
(A,B)_K=(A^*,B^*)_{J_K}= \int_M\tr (A^*\circ B^*)d\mu(g_0)=
$$
$$
=4\int_M\mbox{tr}((1-K^2)^{-1}A(1-K^2)^{-1}B)d\mu(g_0).
$$

To the form $\Omega$:
$$
\Omega_K(A,B)=\Omega_J(A^*,B^*)=\left(A^*J,B^*\right)_J =
\left((AJ_0)^*,B^*\right)_J =
$$
$$
=4\int_M \tr ((AJ_0)^*B^*)d\mu = 4\int_M \tr
\left((1+K)^{-1}AJ_0(1-K)^{-1}(1+K)^{-1}B(1-K)^{-1}\right)d\mu=
$$
$$
=4\int_M \tr \left((1-K^2)^{-1}AJ_0(1-K^2)^{-1}B\right)d\mu.
$$

\vspace{3mm}

Show, that almost complex structure ${\bf J}$ on $\mathcal{A}$ is
integrable. Notice for this, that model space
$\mathrm{End}_{J_0}(TM)$, using for parameterization of the space
$\mathcal{A}$ in neiborhood of point $J_0$ has complex structure
too:
$$
\mathrm{End}_{J_0}(TM)\to \mathrm{End}_{J_0}(TM), \qquad A\to
A\circ J_0.
$$
It seems, that these structures coincide.

\vspace{3mm}

{\bf Theorem 1.} {\it Almost complex structure ${\bf J}$ on
manifold $\mathcal{A}$ is integrable. Corresponding complex
structure coincides with complex structure on $\mathcal{A}$,
obtained by parameterization $\Phi$.}

\vspace{3mm}

{\bf Proof.} Let $A\in \mathrm{End}_{J_0}(TM)$. Then
$d\Phi^{-1}_K(AJ_0)=2J_0(1-K)^{-1}AJ_0(1-K)^{-1}=
2J_0(1-K)^{-1}A(1-K)^{-1}(1-K)J_0(1-K)^{-1}=d\Phi_K(A)\circ J_K$.

\vspace{3mm}

Notice, that weak Pseudo-Riemannian structure on $\mathcal{A}$ is Hermitian
with respect to complex structure ${\bf J}$.
Really, if $A,B\in T_J\mathcal{A}$ are any tangent
elements, then they anticommutating with $J$ and we obtain
$$
\left({\bf J}(A),{\bf J}(B)\right)_J= \int_M \tr (AJBJ)d\mu=
\int_M \tr (AB)d\mu=\left(A,B\right)_J.
$$
Fundamental form of Hermitian weak Pseudo-Riemannian structure
$\left(A,B\right)_J$ on $\mathcal{A}$ coincides with form
$\Omega$ defined as above.

\vspace{3mm}

{\bf Theorem 2.} {\it Fundamental form $\Omega_J$ on
$\mathcal{A}$ is closed.}

\vspace{3mm}

{\bf Proof.} Show, that exterior differential $\Omega_J$ is equal
to zero at any point $J_0\in\mathcal{A}$, $d\Omega_{J_0}=0$. Use
coordinates on $\mathcal{A}$: $J=J_0(1+K)(1-K)^{-1}$. $K=0$
corresponds to the $J_0$. It is enough to show, that
$d\Omega_{K}=0$ when $K=0$. Use standard formula for exterior
differential:
$$
d\Omega(A_0,A_1,A_2) =
A_0\Omega(A_1,A_2)-A_1\Omega(A_0,A_2)+A_2\Omega(A_0,A_1)-
$$
$$
-\Omega([A_0,A_1],A_2)+ \Omega([A_0,A_2],A_1) +
\Omega([A_1,A_2],A_0).
$$
Let $A_0,A_1,A_2$ are constant vector fields on the space
$\mathrm{End}_{J_0}(TM)$. Then all the Lie brackets are equal to
zero. Show that other addends $A_i\Omega(A_j,A_k)$ are equal to
zero too.

Field $A^*$ on $\mathcal{A}$ corresponds to $A\in \mathrm{End}_{J_0}(TM)$
by following formula:
$$
A \mapsto A^*=2J(1+K)^{-1}A(1-K)^{-1}.
$$
Then ${\bf J}(A^*)=(AJ_0)^*=2(1+K)^{-1}AJ_0(1-K)^{-1}$. We obtain
expression of $\Omega$ in coordinate map:
$$
\Omega_K(A,B)=\Omega_J(A^*,B^*)=\left(A^*J,B^*\right)_J =
\left((AJ_0)^*,B^*\right)_J =
$$
$$
=4\int_M \tr ((AJ_0)^*B^*)d\mu = 4\int_M \tr
\left((1+K)^{-1}AJ_0(1-K)^{-1}(1+K)^{-1}B(1-K)^{-1}\right)d\mu=
$$
$$
=4\int_M \tr \left((1-K^2)^{-1}AJ_0(1-K^2)^{-1}B\right)d\mu.
$$
Let $A=A_1,\ B=A_2$ are constant operators (i.e. they are not
depend on $K$). By linearity of integral and trace it is enough to
differentiate expression $(1-K^2)^{-1}A_1J_0(1-K^2)^{-1}A_2$ with
respect to $K$ for finding derivative of $A_0\Omega(A_1,A_2)$. One
can think that $K_t=tA_0$. As
$\left.\frac{d}{dt}\right|_{t=0}(1-K_t^2)^{-1} =
\left.\frac{d}{dt}\right|_{t=0}(1-t^2A_0^2)^{-1}=0$, then
$A_0\Omega(A_1,A_2)=\left.\frac{d}{dt}\right|_{t=0}\Omega_{K_t}(A_1,A_2)=0$.
The theorem is proved.

\vspace{3mm}

{\bf Remark.} Manifold $\mathcal{A}$ isn't Kahlerian, because the inner
product on $\mathcal{A}$ is not positively defined.

\vspace{3mm}

Study the matter of curvature of the space $\mathcal{A}$.

\vspace{3mm}

{\bf Theorem  3.} {\it Space $\mathcal{A}$ has following
geometrical characteristics (in local map $\Phi$).

1) Inner product is given by formula:
$$
(A,B)_K=4\int_M \tr \left((1-K^2)^{-1}A(1-K^2)^{-1}B \right)d\mu,
$$
where $A,B \in \mathrm{End}_{J_0}(TM)$.

2) Covariant derivative of vector fields given by
(constant) operators $A$ and $B$:
$$
\nabla _A B =AK(1-K^2)^{-1}B+BK(1-K^2)^{-1}A.
$$

3) Curvature tensor:
$$
R(A,B)C=-(1-K^2)\left
[\left[(1-K^2)^{-1}A,(1-K^2)^{-1}B\right],(1-K^2)^{-1}C\right],
$$
where $A,B,C\in \mathrm{End}_{J_0}(TM)$.

4) Geodesics, which have beginning at point $J_0$ in directions $A\in
\mathrm{End}_{J_0}(TM)$ are the curves $K(t)$ on
domain $\mathcal{V}(J_0)$:
$$
K(t)=\tanh(t/2 A)=\left( e^{t/2 A}+e^{-t/2 A}\right)^{-1}\left(
e^{t/2 A}-e^{-t/2 A}\right).
$$}

{\bf Proof.} Property 1) is already shown. 2) Covariant
derivative is calculated on six-term formula
\cite{GR}. As operators $A,B,C$ are constant this
formula becames:
$$
(\nabla _A B,C)_K= \frac 12\left(A(B,C)_K + B(A,C)_K -
C(A,B)_K\right).
$$
In calculations we use that integral and trace are linear and
we use the formula:
$$
\left((1-K_t^2)^{-1}\right)'=(1-K^2)^{-1}(AK+KA)(1-K^2)^{-1},
$$
where $K_t=K+tA$ is variation of operator $K$ in direction $A$. Then,
for example,
$$
A(B,C)_K=4\left.\frac{d}{dt}\right|_{t=0}\left( \int_M\tr
((1-K_t^2)^{-1}B(1-K_t^2)^{-1}C)d\mu\right)=
$$
$$
=4\int_M\mbox{tr}((1-K^2)^{-1}((AK+KA)(1-K^2)^{-1}B+
B(1-K^2)^{-1}(AK+KA))(1-K^2)^{-1}Cd\mu
$$
Applying six-term formula we obtain:
$$
(\nabla_AB,C)_K=\frac12\int_M\mbox{tr}((1-K^2)^{-1}(AK(1-K^2)^{-1}B+
BK(1-K^2)^{-1}A+BK(1-K^2)^{-1}A+
$$
$$
+AK(1-K^2)^{-1}B))(1-K^2)^{-1}C)\ d\mu.
$$
Therefore
$$
\nabla_AB=AK(1-K^2)^{-1}B+BK(1-K^2)^{-1}A
$$
for constant vector fields.

3) Curvature tensor is found by formula
$R(A,B)C=\nabla_A\nabla_BC-\nabla_B\nabla_AC+\nabla_{[A,B]}C$. As
$A$ and $B$ are constant vector fields, then $[A,B]=0$. Covariant
derivative $\nabla_BC=BK(1-K^2)^{-1}C+CK(1-K^2)^{-1}B$ depends on
point $K$, then
$\nabla_A\nabla_BC=d_A\nabla_BC+\Gamma(A,\nabla_BC)$.
$$
d_A\nabla_BC=\left.\frac{d}{dt}\right|_{t=0}(BK_t(1-K_t^2)^{-1}C+
CK_t(1-K_t^2)^{-1}B)= BA(1-K^2)^{-1}C+
$$
$$
+BK(1-K^2)^{-1}(AK+KA)(1-K^2)^{-1}C+ CA(1-K^2)^{-1}B+
$$
$$
+CK(1-K^2)^{-1}(AK+KA)(1-K^2)^{-1}B.
$$
We obtain
$$
R(A,B)C=(1-K^2)\left[(1-K^2)^{-1}C,\
[(1-K^2)^{-1}A,(1-K^2)^{-1}B]\right].
$$

4) Direct verification shows that curve
$K(t)=\tanh(t/2 A)$  satisfies to equation $K''=-\Gamma(K',K')$
and, that is why, it is geodesic on coordinate domain
$\mathcal{V}(J_0)$ of space $\mathcal{A}$.

{\bf Remark.} The theorem gives expression of geodesics
$K(t)=\tanh(t/2 A)$ in local map in neiborhood of $J_0$. On the
space $\mathcal{A}$ geodesics are following:
$$
J_t=J_0e^{tA},
$$
where $AJ_0=-J_0A$.

\vspace{3mm}

{\bf 3. Submanifold $\mathcal{AM}$ of associated almost complex
structures.} Suppose, that manifold $M$ is symplectic. It means,
that closed non-degenerated 2-form $\omega $ of class $C^{\infty}$
is given on $M$. The manifold $M$ has even dimension, $\dim M=2n$.
It is natural to research metrics and almost complex structures
which are associated with symplectic form $\omega $.

\vspace{3mm}

{\bf Definition 2.} {\it Almost complex structures $J$ on $M$ is
called \textit{positive associated} with symplectic form $\omega
$, if for any vector fields $X,Y$ on $M$ conditions:

1) $\omega (JX,JY)=\omega (X,Y)$,

2) $\omega (X,JX)>0$, if $X\neq 0$.\\
are hold.}

\vspace{3mm}

{\bf Definition 3.} {\it Every positive associated almost complex
structures $J$ defines Riemannian metric $g$ on $M$ by equality
$$
g(X,Y)=\omega (X,JY),
$$
which is also called \textit{associated}.}

\vspace{3mm}

Associated metric $g$ has following properties:

1) $g(JX,JY)=g(X,Y)$,

2) $g(JX,Y)=\omega (X,Y)$.

\vspace{3mm}

Let $\mathcal{A}_\omega $ is space of all smooth positive
associated almost complex structures on symplectic manifold
$M^{2n},\omega $ and $\mathcal{AM}$ is space of all smooth
associated metrics.

Given spaces $\mathcal{A}_\omega $ and $\mathcal{AM}$
are spaces of smooth sections of corresponding bundles over $M$.
That is why they are infinitely dimensional smooth $ILH$-manifolds.

Let $J\in \mathcal{A}_\omega $ and $g$ is corresponding (unique)
associated Riemannian structure. It is easy to see \cite{Sm}, that
tangent space $T_J\mathcal{A}_\omega$ to manifold of positive
associated almost complex structures consists of symmetric
endomorphisms $P:TM\to TM$, anticommutating with $J$
$$
T_J\mathcal{A}_\omega = \{P\in \mathrm{End}(TM);\quad PJ=-JP,\ \
g(PX,Y) = g(X,PY)\}.
$$
Such space we will denote $\mathrm{End}_{SJ}(TM)$. Therefore
$T_J\mathcal{A}_\omega=\mathrm{End}_{SJ}(TM)$.

Pseudo-Riemannian structure on $\mathcal{A}$ is positive defined
on the submanifold $\mathcal{A}_\omega$. So $\mathcal{A}_\omega$
is weak Riemannian manifold. Manifold $\mathcal{A}_\omega$ has
\cite{Sm} the same properties as space $\mathcal{A}$ (theorems 1
-- 3).

\vspace{3mm}

{\bf Theorem 4.} {\it Space $\mathcal{A}_\omega$ is totally geodesics
submanifold in $\mathcal{A}$.}

\vspace{3mm}

{\bf Proof} follows from coinciding of covariant derivatives on
both spaces $\mathcal{A}$ and $\mathcal{A}_\omega$
\cite{Sm}.

\vspace{3mm}

{\bf Corollary.} Manifold $\mathcal{A}_\omega$ is Kahlerian.

\vspace{3mm}

{\bf 4. Submanifold $\mathcal{AO}$ of orthogonal almost complex
structures.} Consider matter about difference between positive
associated almost complex structures and other ones.

Let $J_0$ is some positive associated almost complex structures
and $g_0$ is corresponding metric. Space $\mathcal{A}$ of all
almost complex structures is parameterized by endomorphisms $K:
TM\to TM$, anticommutating with $J_0$. As $J_0^T = -J_0$, then we
have
$$
J_0K^T = -K^TJ_0
$$
from equality $KJ_0=-J_0K$.
So operator $K$ splits in sum $K=P+L$ of symmetrical $K$ and antisymmetrical
$L$ endomorphisms, each of which anticommutates with $J_0$,
$$
\mathrm{End}_{J_0}(TM)=\mathrm{End}_{SJ_0}(TM)\oplus\mathrm{End}_{AJ_0}(TM).
$$

At exponential parameterization of space $\mathcal{A}$
$$
E:\ \mathrm{End}_{J_0}(TM) \longrightarrow \mathcal{A},\quad
K\mapsto J = J_0e^K,
$$
subspace $\mathrm{End}_{SJ_0}(TM)$ of symmetrical endomorphisms
parameterizes associated almost complex structures, and subspace
$\mathrm{End}_{AJ_0}(TM)$ of antisymmetrical endomorphisms using
for parameterization of remain, non-associated almost complex
structures.

Thus, submanifold, which is transversal to
$\mathcal{A}_\omega$ is parameterized by mapping
$$
E_A:\ \mathrm{End}_{AJ_0}(TM) \longrightarrow \mathcal{A},\quad
L\mapsto J = J_0e^L.
$$
As endomorphism $L$ is antisymmetrical, then $e^L$ is orthogonal
transformation anticommutating with $J_0$.

So submanifold, which is transversal to $\mathcal{A}_\omega$ in
neiborhood of $J_0$ formed by almost complex structures $J$, such
that $J=J_0O$, where $O$ is orthogonal transformation,
anticommutating with $J_0$. Such almost complex structures are
orthogonal and they form submanifold in space
$\mathcal{A}_\omega$.

Let $J_0$ is some positive associated almost complex structure and
$g_0$ is corresponding metric.

Consider the set $\mathcal{AO}(M)$ $g_0$ of orthogonal almost
complex structure $J$:
$$
\mathcal{AO}(M)=\{J\in \mathcal{A}:\ g_0(JX,JY)=g_0(X,Y),\ J\
\mbox{ gives the same orientation as }\ J_0 \}
$$

This set is $ILH$-submanifold of the space $\mathcal{A}$.

It is easy to see that tangent space $T_J\mathcal{AO}(M)$
consists of antisymmetrical endomorphisms $K:TM\to TM$,
anticommutating with $J$
$$
T_J\mathcal{AO})(M) = \{K\in \mathrm{End}(TM);\quad KJ=-JK,\ \
g_0(KX,Y) = - g_0(X,KY)\}.
$$
We will denote such space as $\mathrm{End}_{AJ}(TM)$. Thus
$T_J\mathcal{AO}(M)=\mathrm{End}_{AJ}(TM)$.

Almost complex structure ${\bf J}$ on $\mathcal{A}$ lefts
submanifold $T_J\mathcal{AO}(M)\subset T_J\mathcal{A}$ invariant.
Therefore $\mathcal{AO}(M)$ is complex submanifold in
$\mathcal{A}$. Pseudo-Riemannian structure on $\mathcal{A}$ is
negative defined on $\mathcal{AO}(M)$. So $\mathcal{AO}(M)$ can be
viewed as weak Riemannian manifold.

It is easy to see, that all computations in the proof of theorem 3
can be repeated in this case. Therefore manifold
$\mathcal{AO}(M)$ has the same properties as the space
$\mathcal{A}$ (theorems 1 -- 2).

\vspace{3mm}

{\bf Theorem 5.} {\it Submanifold $\mathcal{AO}(M)$ is
totally geodesic submanifold in $\mathcal{A}$.}

\vspace{3mm}

{\bf Proof} follows from coinciding of covariant derivatives on
spaces $\mathcal{AO}(M)$ and $\mathcal{A}$,
this checked by direct calculations.

\vspace{3mm}

{\bf Corollary.} Manifold $\mathcal{AO}(M)$ is Kahlerian.

\end{document}